\documentclass[11pt]{article}
\usepackage{latexsym}

\newcommand{\ms}[1]{\mbox{\hspace{#1cm}}}
\newcommand{\fr}[1]{\stackrel{#1}{\longrightarrow}}
\newcommand{\A}{\stackrel{2}{\wedge}}
\newcommand{\x}{\otimes}
\newcommand{\p}{\oplus}
\newcommand{\fd}{\rightarrow}

\newcommand{\ac}[1]{\tilde{#1}}
\newcommand{\Om}{\Omega^1}

\newcommand{\Ho}{H^{0}}

\newcommand{\ch}{\vee}

\newtheorem{theorem}{Theorem}[section]

\newtheorem{cor}{Corollary}[section]

\newtheorem{pr}{Proposition}[section] 
\newenvironment{prop}{\begin{pr} \mbox{}\\[2mm]}{\end{pr}}

\newtheorem{lemma}{Lemma}[section]
\newenvironment{lem}{\begin{lemma} \mbox{}\\[2mm]}{\end{lemma}}
\newtheorem{df}{Definition}[section]
\newenvironment{defin}{\begin{df} \mbox{}\\[2mm] \slshape }{$\Box$\end{df}}
\newtheorem{os}{Remark}[section]
\newenvironment{oss}{\noindent \begin{os} \mbox{}\\}{\end{os}}
\newenvironment{proof}{\noindent{\em proof.}}{$\Box$}
 
\newcounter{lett}

\newcounter{splist}

\newcounter{dslist}

\newcommand{\PP}{\mbox{{I \hspace{-,18cm}P}}}
\newcommand{\HH}{\mbox{{I \hspace{-.18cm}H}}}
\newcommand{\TPP}{\mbox{{I \hspace{-,18cm}P}}}
\newcommand{\TTPP}{\mbox{{I \hspace{-,32cm}P}}}
\newcommand{\KK}{\mbox{{\itshape I \hspace{-.38cm} K}}}
\newcommand{\ol}[1]{\overline{#1}}

\newcommand{\hs}[1]{\hspace{#1}}

\newcommand{\refeq}[1]{(\mbox{\ref{#1})}}

\begin{document}
\begin{center}
{\LARGE \bf
 
 The Tangent  Space at a Special Symplectic Instanton 
   Bundle on $\TTPP^{2n+1}$ } \\
 CARLA DIONISI \\
   \end{center}
 
 {\small \it {\bf  ABSTRACT }
  :Let $MI_{\mbox{Simp},\PP^{2n+1}}(k)$ be the moduli space of stable 
symplectic instanton bundles on $\PP^{2n+1}$ with second Chern class $c_2=k$
(it is a closed subscheme of the moduli space $MI_{\PP^{2n+1}}(k)$) \\
We prove that the dimension
of its Zariski tangent space at a special (symplectic) instanton bundle is \
 $2k(5n-1)+4n^2-10n+3 \quad , k\geq 2$. \ 
 It follows that special symplectic instanton bundles are smooth points for 
 $ k \leq 3 $}

\section*{Introduction}

Symplectic instanton bundles on $\TPP^{2n+1}$ are holomorphic bundles of rank 2n
(see \cite{A} ,\cite{MS}  and \cite{OS}) that correspond to the self-dual solutions of 
Yang-Mills equations on $\TPP^n(\HH)$. \
They are given by some monads (see section 2
for precise definitions) and their only topological invariant is  $c_2=k$.\\
At present the dimension of their moduli space $MI_{\mbox{Simp},\PP^{2n+1}}(k)$
is not known except in the cases n=1, where the dimension is $8k-3$
(see \cite{HN}), and in 
few other cases corresponding to small values of k.\\
$MI_{\mbox{Simp},\PP^{2n+1}}(k)$ is a closed subscheme of 
$MI_{\PP^{2n+1}}(k)$ and this last scheme parametrizes stable instanton
bundles with structural group $GL(2n)$.\\
The class of special instanton bundles was introduced in \cite{ST}. \\
Let $E \in MI_{\PP^{2n+1}}(k)$ be a special symplectic instanton bundle. The tangent
dimension $h^1(End(E))$ was computed in \cite{OT} and it is equal to       
$4(3n-1)k+(2n-5)(2n-1)$ .\\
The Zariski tangent space of $MI_{\mbox{Simp},\PP^{2n+1}}(k)$
at $E$ is isomorphic to $H^1(S^2E)$
and in this paper we prove that
 
\begin{equation}
\label{zero}
h^2(S^2E)=\left(\begin{array}{c} k-2 \\ 2 \end{array} \right) \cdot 
 \left(\begin{array}{c} 2n-1 \\ 2 \end{array} \right) \qquad   \forall \ k \geq 2
 \end{equation}
 
By the Hirzebruch-Riemann-Roch formula , since $h^0(S^2E)=0$ and $h^i(S^2E)=0 \ 
\forall\ i\geq3$, it follows that:
\begin{center}
$\chi(S^2E)=h^2(S^2E)-h^1(S^2E)=2n^2+n+\frac{1}{2}\left[ k^2\left( \begin{array}{c} 2n-1 \\ 2 \end{array} \right) -k(10n^2-5n-1) \right] $
\end{center}
and
 \begin{theorem} \label{acca1}
 Let $E$ be a special symplectic instanton bundle.Then
 $$h^1(S^2E)=2k(5n-1)+4n^2-10n+3 \qquad  ,k\geq 2 $$
 \end{theorem}
 (for $n=1$ it is well known that \ $ h^1(S^2E)=8k-3$ ).\\
 Now, since by the Kuranishi map $H^2(S^2E)$ is the space 
of obstructions to the smoothness at $E$ of $MI_{\mbox{Simp},\PP^{2n+1}}(k)$,
we obtain
\begin{cor}
$\forall k\geq 2$
 the dimension of any irreducible component of
 $MI_{\mbox{Simp},\PP^{2n+1}}(k)$,
containing a special symplectic instanton bundle is bounded by the value \\
$2k(5n-1)+4n^2-10n+3 \quad $ ( linear in k )

\end{cor}

\begin{cor}
 $\forall n$ \  $MI_{\mbox{Simp},\PP^{2n+1}}(3)$ is smooth at a special instanton
 bundle $E$,and the dimension of any irreducible component containing $E$ is
 $ 4n^2+20n-3 $.
 \end{cor}
 The main remark of this paper is that it is easier to compute $H^2(S^2E)$
 and $H^2(\A E)$ together as SL(2)-modules (although this second cohomology 
 space has a  geometrical meaning only for orthogonal bundles) than to compute
  $H^2(S^2E)$ alone.
 
 \section{Preliminaries}
 Throughout this paper $\KK$ denotes an algebraically closed field of characteristic
 zero.
   $U$ denotes a 2-dimensional $\KK$ vector space \((U=<s,t>)\), \ $S_n=S^nU$ 
   its n-th symmetric power   $\ms{0.2}(S_n=<s^n, s^{n-1}t,\ldots, t^n>$) ,
  $V_n=U\x S_n $ $\ms{0.2}(V_n=<s\x s^n, s\x s^{n-1}t,\ldots s\x t, \ldots t\x t^n>)$ \ 
  and $\PP^{2n+1}=\PP(V_n)$.

  \begin{defin}
  \label{def 1.1}
  A vector bundle $E$ on $\PP^{2n+1}$ of rank $2n$ is called an {\bf instanton bundle of
   quantum number $k$} if:
\begin{itemize}
\pagebreak

\item $E$ has Chern polinomial $c_t(E) \ = \ (1-t^2)^{-k}$;
\item $E(q)$ has  natural cohomology in the range $-(2n+1) \ \leq \ q \ \leq 0$,
that is $H^i(E(q)) \ \neq \ 0$ for at most one $i = i(q)$.
\end{itemize}
\end{defin}

By \cite{OS},\cite{AO1},the Definition ~\ref{def 1.1} is
equivalent to : \\
i)$E$ is the cohomology bundle of a monad:
\[
0 \fd O(-1)^k \fd \Omega^1(1)^k \fd O^{2n(k-1)} \fd 0
\]
 or ii) $E$ is the cohomology bundle of a monad:
\[
 0 \fd O(-1)^k \fr{A} O^{2n+2k} \fr{B^t} O(1)^k \fd 0
\]
(where, after we have fixed a coordinate system, A and B can be identified 
with matrices in the space 
$Mat(k,2n+2k,S_1)$)

\begin{defin}
 An instanton bundle $E$
  is called {\bf symplectic} if there is 
 an isomorphism $\varphi:E \fd E^ {\ch}$ satisfying
  $\varphi = -\varphi^{\ch}$.
 \end{defin}

\begin{defin}
 An instanton bundle is called {\bf special} if it arises from a monad 
 where the morfism $B^t$
 is defined in some system of homogeneous coordinates $x_0, \cdots x_n, y_0 \cdots y_n$ 
 on $\PP^{2n+1}$ by the trasposed of the matrix:
  \[
B = \left(
\begin{array}{cccccccccccc}
x_0 & \cdots & x_n & 0 & \cdots & 0 & y_0 & \cdots & y_n & 0 & \cdots & 0 \\
0 & x_0 & \cdots & x_n & 0 & \cdots & 0 & y_0 & \cdots & y_n & 0 & \cdots \\
\cdots & \cdots & \cdots & \cdots & \cdots & \cdots & \cdots & \cdots & \cdots & \cdots & \cdots & \cdots\\
\cdots & 0 & x_0 & \cdots & x_n & 0 & \cdots & 0 & y_0 & \cdots & y_n & 0 \\
0 & \cdots & 0 & x_0 & \cdots & x_n & 0 & \cdots & 0 & y_0 & \cdots & y_n
\end{array}
\right)
\] \\
 \end{defin}
The following lemma is well known (and easy to prove)

\pagebreak

\begin{lemma}
$$\Ho (O(1)) \cong V^{\ch}$$ 
$$\Ho(\Om(2)) \cong \A V^{\ch}$$
$$H^i (\PP^n, S^2 \Om(1)) = \left \{ \begin{array}{ll}
                        0 & \mbox{se } i\neq 1 \\
                        \A V^{\ch} & \mbox{se } i= 1
                        \end{array} \right.$$
  \end{lemma}
\section{Existence of a special symplectic instanton bundle}
There is a natural exact sequence of GL(U)-equivariant maps for any $k , n\geq 1$
(Clebsch-Gordan sequence):
 \begin{equation}
\label{uno_3}
 0 \fd \A U \x S_{k-1} \x V_{n-1} \fr{\beta} S_k \x V_n \fr{\mu} V_{k+n} \fd 0
\end{equation}
where $\mu$ is the multiplication map and $\beta$ is defined by
$(s \wedge t) \x f \x g \fd (sf \x tg - tf \x sg)$

 We can define (see \cite{OT}) the morphism \\
 \centerline{
 $\ac{b}: S_{k-1}^{\ch} \x \Om(1) \fd \A U^{\ch}  \x S_{k-2}^{\ch} \x V_{n-1}^{\ch} \x O$
 }
 
 and it is induced the complex 
 \begin{equation}
\label{cinque_3}
A \x O(-1) \fr{\ac{a}} S^{\ch}_{k-1} \x \Om(1) \fr{\ac{b}} \A U^{\ch} \x S^{\ch}_{k-2} \x V_{n-1}^{\ch} \x O
\end{equation}
  where A is a k-dimensional subspace of $S^{\ch}_{2n+k-1} \x \A U^{\ch}$ 
  such that ~\refeq{cinque_3} is a monad and the cohomology bundle $E$ is a special symplectic
  instanton bundle.
It was proved in \cite{OT} that 
$$ H^2(EndE)\cong Ker(\Phi^{\ch})^{\ch}$$
where
 \centerline{
  $\Phi^{\ch}: S_{k-2}^{\x 2} \x V_{n-1}^{\x 2} \fd S_{k-1}^{\x 2} \x \A V_{n}$
  }
and there is an isomorphism of 
SL(2)-representations 
$$ 
\varepsilon : S^{\ch}_{k-3} \x S^{\ch}_{k-3} \x S^2V^{\ch}_{n-2} \fd Ker(\Phi^{\ch})$$

\section{How to identify $H^2(S^2E)$ and  $H^2(\A E)$} 
 \begin{prop}
 Let $E$ be special symplectic instanton bundle ,cohomology of monad
  ~\refeq{cinque_3} and  $N=Ker \ac{b}$.Then   
 \begin{description}
\item[(i)] $H^2 (S^2E) \cong H^2(S^2N) $
\item[(ii)] $H^2 (\A E) \cong H^2(\A N) $
\end{description}
\end{prop}

\begin{proof}   \\
We denote $B:= S_{k-1}^{\ch}$ \qquad and \qquad $C:=\A U^{\ch} \x S_{k-2}^{\ch}  \x V_{n-1}^{\ch} $ \\
The result follows from the two exact sequences given by 
 monad ~\refeq{cinque_3} :
\begin{equation}
\label{cinque_4}
 0 \fd N \fd B \x \Om (1) \fd C \x O \fd 0
\end{equation}

\begin{equation}
0 \fd A \x O(-1) \fd N \fd E \fd 0
\end{equation}
In fact, by performing the second symetric and alternating power
of sequence ~\refeq{cinque_4}, we have 

\centerline{
\begin{minipage}{2in}
\begin{tabbing}
$ 0 \fd S^2 N \fd $\=$\tilde{A}$\=$\fd B \x C \x \Om(1) \fd \A C \x O \fd 0$\\
                   \>            \>$\searrow$\=\hspace{0.5 cm}$\nearrow$\=\\
                   \>            \>          \> $ M^1$\\ 
                   \>            \>$\nearrow$\hspace{0.5 cm}$\searrow$\\
                   \>       0    \>          \>                     \> 0
\end{tabbing}
\end{minipage}
}
     
\begin{equation}
\label{DUE}
\end{equation}
where $\tilde{A}:= S^2(B \x \Om(1)) = (S^2 B \x S^2(\Om(1))) \p ( \A B \x \Omega^2(2))$ \\
and \\
\centerline{
\begin{minipage}{2in}
\begin{tabbing}
$ 0 \fd \A N \fd $\=$\overline{A} $\=$\ \fd B \x C \x \Om(1) \fd O \x S^2 C \fd 0$\\
                   \>            \>$\searrow$\=\hspace{0.5 cm}$\nearrow$\=\\
                   \>            \>          \> $ M$\\ 
                   \>            \>$\nearrow$\hspace{0.5 cm}$\searrow$\\
                   \>       0    \>          \>                     \> 0
\end{tabbing}
\end{minipage}
}     
\begin{equation}
\label{TRE}
\end{equation}
where $\overline{A}:= \A(B \x \Om(1)) = (\A B \x S^2(\Om(1))) \p ( S^2 B \x \Omega^2(2))$ \\
\end{proof}
\subsection{Identifying $ H^2(S^2 N)$ and $H^2 ( \A N) $}
i)
\label{par43}
 Diagram ~\refeq{DUE} gives the following two exact sequences:\\

\begin{minipage}{5 in}
\begin{equation}\label{quattro4} 
O \fd H^0(M^1) \fd H^1(S^2N) \fd H^1(\ac{A}) \fd H^1(M^1) \fd H^2(S^2(N)) \fd H^2(\ac{A}) \fd \cdots
\end{equation}
\begin{tabbing}
$O \fd H^0(M^1) \fd B \x C \x H^0$\=$(\Om(1))\fd \A C\fd H^1(M^1) \fd B \x C \x H^1$\=$(\Om(1))\fd \cdots$\\
                     \>$\parallel$\> $\parallel$\\ 
                     \>$0 $       \> $0$
\end{tabbing}
\begin{equation}\label{cinque4} 
\end{equation} 
\end{minipage}
 
Sequence \refeq{cinque4} implies:

$H^0(M^1) = 0$   \quad and \quad $H^1(M^1) \cong \A C$

Then, by using the two formulas:

$H^1(\tilde{A}) = (S^2 B \x H^1(S^2 \Om(1))) \p ( \A B \x H^1(\Omega^2(2)) 
=S^2 B \x \A V^{\ch}$

and:\

$H^2(\tilde{A})=  (S^2 B \x H^2(S^2 \Om(1))) \p ( \A B \x H^2(\Omega^2(2)) =0$

sequence ~\refeq{quattro4} becomes:

\[0 \fd H^1(S^2N) \fd H^1(\ac{A}) \fd H^1(M^1) \fd H^2(S^2(N)) \fd 0\]

i.e.\[0 \fd H^1(S^2 N) \fd S^2 B \x \A V^{\ch} \fr{\ac{\Phi}} \A C \fd H^2(S^2N) \fd 0 \]

\[ \Longrightarrow \qquad H^2 ( S^2 N ) \cong \mbox{Coker}(\ac{\Phi}) = (\mbox{Ker}(\ac{\Phi}^{\ch}))^{\ch} \]

Then:
\[H^2 ( S^2 N) ^{\ch} = \mbox{Ker} \left[ \A (S_{k-2} \x V_{n-1} ) \fr{\ac{\Phi}^{\ch}} S^2(S_{k-1}) \x \A V_n\right]\]

ii)                                 
Diagram  ~\refeq{TRE} gives the following two exact sequences:

\begin{minipage}{5in}
 \begin{equation}\label{sei4} 
O \fd H^0(M) \fd H^1(\A N) \fd H^1(\overline{A}) \fd H^1(M) \fd H^2(\A N) \fd H^2(\overline{A}) \fd \cdots
\end{equation}

\begin{tabbing}
$O \fd H^0(M) \fd B \x C \x H^0$\=$(\Om(1))\fd S^2 C $\=$\x H^0(O) \fd H^1(M) \fd 0 \fd \cdots$\\
                     \>$\parallel$\> $\parallel$\\ 
                     \>$0 $       \> $S^2 C$
\end{tabbing}
\begin{equation}\label{sette4} 
\end{equation} 
\end{minipage}\\

and, from sequence ~\refeq{sette4},  we get \\
$H^0(M) = 0$ \quad and \quad $H^1(M) \simeq S^2 C$ \\
Then, since :\\
$H^1(\overline{A})=(H^1(S^2(\Om(1))\x \A B) \p (S^2 B \x H^1(\Omega^2(2)))=
       \A B \x \A V^{\ch}$ \\
and \qquad $H^2(\overline A) = 0$\\
sequence  ~\refeq{sei4} becomes :

\begin{minipage}{5 in}
\begin{tabbing} 
$O \fd H^0$\=$(M) \fd H^1(\A N)$\=$ \fd H^1(\overline{A}) \fd H^1(M) \fd H^2(\A N) \fd 0$\\
\>$\parallel$\\
\>0
\end{tabbing}
\end{minipage}

i.e. \qquad
$0 \fd H^1(\A N) \fd \A B \x \A V^{\ch} \fr{\overline{\Phi}} S^2 C \fd H^2(\A N ) \fd 0 $ \\
\[
\Longrightarrow \qquad  H^2 ( \A N ) \cong \mbox{Coker}(\overline{\Phi}) = (\mbox{Ker}(\overline{\Phi}^{\ch}))^{\ch} \]
Then we obtain :
\[(H^2(\A N))^{\ch} = \mbox{Ker} \left[ S^2(S_{k-2}\x V_{n-1})\fr{\overline{\Phi}^{\ch}}\A S_{k-1}\x \A V_n \right]\]

\subsection{\bf Identifying $H^2(S^2 E)$}
We have \qquad \qquad  $H^2(S^2E)^{\ch} \cong \mbox{Ker} \ \ac{\Phi}^{\ch}$ \\
     where \qquad
$\ac{\Phi}^{\ch} : \A(S_{k-2}\x V_{n-1}) \fd S^2 S_{k-1} \x \A V_n$ \\
is explicitly given by
\begin{tabbing}
$\ac{\Phi}^{\ch}((g \x v)\wedge (g^1 \x v^1)) =$\=$sg \cdot sg^1 \x (tv \wedge tv^1)-sg\cdot tg^1 \x (tv \wedge sv^1)+$\\
\>$-tg\cdot sg^1 \x (sv \wedge tv^1)+tg\cdot tg^1 \x (sv \wedge sv^1) $ \\
\end{tabbing}
i.e. \qquad $ \ac{\Phi}^{\ch} = \ac{p} \circ (\A\beta)$ \\
where  \qquad
$\beta : \ \ \A U \x S_{k-2} \x V_{n-1} \fd S_{k-1} \x V_n \; \; \; $
is such  that
\[ 
(s \wedge t) \x (g \x v) \mapsto (sg \x tv) - (tg \x sv)
\]
and \begin{tabbing}
$\ac{p} \ : \ $\=$\A($\=$S_{k-1} \x V_n) \fd S^2 S_{k-1} \x \A V_n$ \\
\>\>$\|$ \\
\>$(\A S_{k-1} \x S^2V_n) \p (S^2S_{k-1} \x \A V_n)$\\
\end{tabbing}
is such that
$$
(f \x u) \wedge (f' \x u^1) \mapsto f\cdot f' \x u \wedge u^1 \mbox{.}
$$
Now, we consider the $SL(2)$-equivariant morphism:\\
\[ 
\ac{\varepsilon}^1 \ : \ \A (S_{k-3} \x V_{n-2}) \fd \A(S_{k-2} \x V_{n-1})
\]
where, up to the order of factors, the map
$\ac{\varepsilon}^1 := \beta^1 \wedge \beta^1$ and $\beta^1:S_{k-3} \x V_{n-2} \fd S_{k-2} \x V_{n-1}$ 
is defined as  $\beta$.  Hence, $\ac{\varepsilon}^1 $ is injective.

Finally, we define
\[ 
\ac{\varepsilon} \ : \ \A S_{k-3} \x S^2V_{n-2} \fd \A(S_{k-2} \x V_{n-1})
\]

as $\ac{\varepsilon} = \ac{\varepsilon}^1 \circ \ac{i}$, \  where
\begin{tabbing}
$\ac{i} \ :\ $\=$ \A S_{k-3} \x S^2V_{n-2} \fd \A(S_{k-3} \x V_{n-2}) \ \ $ such that \\
\> $f \wedge f' \x u\cdot u^1 \longmapsto (f \x u) \wedge (f' \x u^1) + (f \x u^1) \wedge (f' \x u)$\\
\end{tabbing}
is an injective map. Then,
also $\ac{\varepsilon}$ is injective.

\begin{lem} 
\label{lemma442}
Im $ \ac{\varepsilon} \subset \mbox{Ker} \ \ac{\Phi}^{\ch}$ 
 \end{lem}
 \begin{proof}
Straightforward computation.
\end{proof}

\subsection {\bf Identifying $H^2(\A E)$}

We have  $$H^2(\A E)^{\ch} \cong \mbox{Ker} \ \ol{\Phi}^{\ch}$$
where  \qquad
$
\ol{\Phi}^{\ch} \ : \ S^2(S_{k-2}\x V_{n-1}) \fd \A S_{k-1} \x \A V_n
$ \qquad
is explicity given by
\begin{tabbing}
$\ol{\Phi}^{\ch}((g \x v)\cdot (g^1 \x v^1)) =$\=$sg \wedge sg^1 \x (tv \wedge tv^1) - sg \wedge tg^1 \x(tv \wedge sv^1)$\\
\>$-tg \wedge sg^1 \x (sv \wedge tv^1)+(tg \wedge tg^1) \x (sv \wedge sv^1)$
\end{tabbing}
i.e. \qquad
 $\ol{\Phi}^{\ch} = \ol{p} \circ(S^2 \beta)$ \qquad where
\pagebreak
\begin{tabbing}
$\ol{p} \ : \ $\=$S^2($\=$S_{k-1}\x V_n) \fd \A S_{k-1} \x \A V_n$\\
\>\>$\|$\\
\>$(\A S_{k-1}\x \A V_n) \p (S^2 S_{k-1} \x S^2 V_n)$
\end{tabbing}
is such that
\[
\ol{p}((f \x u) \cdot (f' \x u^1)) = f \wedge f' \x u \wedge u^1
\]
We consider the $SL(2)$-equivariant morphism:
\[
\ol{\varepsilon}^1 \ :\ S^2(S_{k-3} \x V_{n-2}) \fd S^2 (S_{k-2} \x V_{n-1})
\]
such that:
\begin{tabbing}
$\ol{\varepsilon}^1((f \x u) \cdot (f' \x u^1)) =$\=$(sf \x tu) \cdot(sf' \x tu^1)-(sf \x su) \cdot (tf' \x tu^1)+$\\
\>$-(tf \x tu) \cdot (sf' \x su^1) + (sf \x tu) \cdot (sf' \x tu^1)$
\end{tabbing}
($\ol{\varepsilon}^1 = S^2 \beta^1 $ hence  $ \ol{\varepsilon}^1$ is injective).
Finally, we define
\[
\ol{\varepsilon} \ : \ S^2 S_{k-3} \x S^2 V_{n-2} \fd S^2 (S_{k-2} \x V_{n-1})
\]
as $\ol{\varepsilon} = \ol{\varepsilon}^1 \circ \ol{i}$ \quad where
\begin{tabbing}
$\ol{i}\ : \ $\=$ S^2 S_{k-3} \x S^2V_{n-2} \fd S^2 (S_{k-3} \x V_{n-2})$ \ such that \\
\>$f \cdot f' \x uu^1 \mapsto (f \x u)(f' \x u^1) + (f \x u^1)(f' \x u)$
\end{tabbing}
is an injective map. Then, also $\ol{\varepsilon}$ is injective
\begin{lem}
\label{lemma443}
Im $ \ol{\varepsilon} \subset \mbox{Ker} \ \ol{\Phi}^{\ch}$ 
\end{lem}
\begin{proof}
Straightforward computation.
\end{proof}
\begin{theorem}
For any special symplectic instanton bundle $E$
\[
H^2(S^2E)\ \simeq\ \A (S_{k-3})^{\ch} \x S^2(V_{n-2})^{\ch}\]
 \end{theorem}
 \begin{proof} 
By lemma ~\ref{lemma442} and ~\ref{lemma443} we have the following diagram
with exact rows and columns:

\begin{minipage}{4in}
\begin{tabbing}
$0 \fd H^2(\A$\=$ N)^{\ch} \fd S^2(S_{k-2} \x$\=$V_{n-1}) \fr{\ol{\Phi}^{\ch}} \A S_{k-1} \x$\=$\A V_n \fd H^1(\A$\=$N)^{\ch}\fd 0$ \kill
\> $0$ \> $0$ \> $0$ \> $0$\\
\>$\downarrow$ \> $\downarrow$ \> $\downarrow$\> $\downarrow$\\
$0 \fd H^2(\A N)^{\ch} \fd S^2(S_{k-2} \x V_{n-1}) \fr{\ol{\Phi}^{\ch}} \A S_{k-1} \x \A V_n \fd H^1(\A N)^{\ch}\fd 0$ \\     
\>$\downarrow$ \> $\downarrow$ \> $\downarrow$\> $\downarrow$\\
$0 \fd H^2(N \x N)^{\ch} \fd S_{k-2} ^{\x 2} \x V_{n-1}^{\x 2} \fr{\Phi^{\ch}} S_{k-1}^{\x 2} \x \A V_n \fd H^1(N \x N)^{\ch}\fd 0$ \\     
\>$\downarrow$ \> $\downarrow$ \> $\downarrow$\> $\downarrow$\\
$0 \fd H^2(S^2 N)^{\ch} \fd \A (S_{k-2} \x V_{n-1}) \fr{\ac{\Phi}^{\ch}} S^2 S_{k-1} \x \A V_n \fd H^1(S^2  N)^{\ch}\fd 0$ \\     
\>$\downarrow$ \> $\downarrow$ \> $\downarrow$\> $\downarrow$\\
\> $0$ \> $0$ \> $0$ \> $0$\\
\end{tabbing}
\end{minipage}    \\
It was shown in \cite{OT} that:

$H^2(EndE) \simeq{Ker}\ \Phi^{\ch} = H^2(N \x N)^{\ch} \simeq S_{k-3}^{\x 2} \x S^2 V_{n-2}$

We have proved that there are two injective maps:
\[
\ac{\varepsilon}\ : \ \A (S_{k-3}) \x S^2V_{n-2} \fd \mbox{Ker} \ \ac{\Phi}^{\ch} \simeq H^2(S^2N)^{\ch} \simeq H^2(S^2E)^{\ch}
\]
\[
\ol{\varepsilon}\ : \ S^2(S_{k-3}) \x S^2 V_{n-2} \fd \mbox{Ker} \ \ol{\Phi}^{\ch} \simeq H^2(\A N)^{\ch} \simeq H^2(\A E)^{\ch}
\]

 Then, we can consider the following diagram:                                        
\begin{tabbing}
$0 \fd S^2 S_{k-3} $\=$ \x S^2 V_{n-2} \fd S_{k-3}^{\x 2} \x $\=$ S^2V_{n-2} \fd \A S_{k-3}$\=$\x V_{n-2}$\kill
\> $0$ \> $0$ \> $0$ \\
\> $\downarrow$ \> $\downarrow$ \> $\downarrow$ \\
$0 \fd S^2 S_{k-3} \x S^2 V_{n-2} \fd S_{k-3}^{\x 2} \x  S^2V_{n-2} \fd \A S_{k-3} \x S^2V_{n-2} \fd 0 $\\
\> $\downarrow \ol{\varepsilon}$ \> $\downarrow  \varepsilon$ \> $\downarrow \ac{\varepsilon}$ \\
\ \ \ \ \ \ \ \ \ $0\fd $ \> $H^2(\A E)^{\ch} \fd H^2(EndE)^{\ch} \fd H^2(S^2E)^{\ch} \fd 0 $ \\
\>\> $\downarrow$\\
\>\> $0$\\
\end{tabbing}
and by the {\bf Snake-Lemma} there is the exact sequence :
\begin{tabbing}
$0 \fd$ Ker \ \= $\ol{\varepsilon} \fd $ Ker\ \= $\varepsilon \fd $ Ker \= $\ac{\varepsilon} \fd \mbox{Coker}\ \ol{\varepsilon} \fd \mbox{Coker}\ $\=$ 
\varepsilon \fd \mbox{Coker}\ \ac{\varepsilon} \fd 0$\\
\> $\|$ \> $\|$ \> $\|$ \> $\|$ \\
\> $0$ \> $0$ \> $0$ \> $0$\\
\end{tabbing}
$
\Rightarrow \mbox{Coker} \ \ol{\varepsilon} = 0 \Rightarrow \ol{\varepsilon}\  \mbox{is an 
isomorphism} \ \Rightarrow \ac{\varepsilon}\ \mbox{is an isomorphism}.
$ \\
Thus:
\[H^2(S^2E)^{\ch}\  \cong\ \A (S_{k-3}) \x S^2(V_{n-2})\]
i.e. \qquad \ 
$H^2(S^2E)\ \simeq\ \A (S_{k-3})^{\ch} \x S^2(V_{n-2})^{\ch}$ \quad as we wanted.\\
\end{proof}
\begin{oss} By
this theorem formula ~\ref{zero} and  theorem ~\ref{acca1} are easily
 proved.
\end{oss}

 Carla Dionisi \\
 Dipartimento di Matematica ed Applicazioni "R.Caccioppoli" \\
 Universit\`{a} di Napoli Federico II \\
 via Cintia (loc. Monte S.Angelo) \\
 I-80138 Napoli, Italy \\
 {\em E-mail address}: dionisi@matna2.dma.unina.it

 \end{document}